\documentstyle{article}
\begin{document}

\title{COORDINATE CALCULI ON ASSOCIATIVE ALGEBRAS\footnote{Invited
talk
presented by V.K. Kharchenko at the $XXX_q$ Winter School in
Karpacz.}}
\author{
Andrzej Borowiec \thanks{Supported by the State Research Committee
KBN No 2 P302 023 07.}\\Institute of Theoretical Physics \\
University of Wroc{\l}aw, Poland, borowiec@ift.uni.wroc.pl
\and
Vladislav K. Kharchenko \thanks{Supported by the Russian Fund of
Fundamental
Research No 93-011-16171.}\\ Institute of Mathematics\\
Novosibirsk, Russia, kharchen@math.nsk.su}
\date{}
\maketitle
\begin{abstract}
A new notion of an optimal algebra for a first order free differential
was introduced in \cite{BKO}. Some relevant examples are indicated.
Quadratic identities in the optimal algebras and calculi on
quadratic algebras are studied. Canonical construction of a
quantum de Rham complex for the coordinate differential is proposed.
The relations between calculi and various generalizations
of the Yang--Baxter equation are established.
\end{abstract}

\def\c{\cdot}
\def\ot{\otimes}
\def\rt{\rightarrow}
\def\ld{\ldots}

\def\ba{\begin{array}}
\def\ea{\end{array}}

\def\N{{I\!\! N}}
\def\F{{I\!\! F}}
\def\R{\hat{R}}
\def\I{\hat{I}}
\def\A{\bar{A}}
\def\x{\hat{x}}
\def\A{\hat{A}}
\def\pih{\hat{\pi}}
\def\D{\hat{D}}
\def\d{\hat{d}}
\def\M{\cal M}

\def\al{\alpha}
\def\be{\beta}
\def\de{\delta}
\def\la{\lambda}
\def\La{\Lambda}
\def\ga{\gamma}
\def\Om{\Omega}
\def\om{\omega}

\section{Introduction.}
Quantum spaces are identified with noncommutative algebras.
This idea generalizes the concept of supersymmetry and opens
new possibilities for quantization of classical systems.
New objects like quantum groups and quantum spaces are usually
considered
as (multi)parametric noncommutative (=quantum) deformations of the
corresponding classical (=commutative) objects.

Differential calculi on quantum spaces have been
elaborated by Pusz and Woro\-no\-wicz \cite{PW,P} and by
Wess and Zumino \cite{WZ}.
Bicovariant differential calculi on  quantum groups were presented by
Woronowicz \cite{W}. Woronowicz found that
construction of first order calculus, a bimodule of one forms
is not functorial. There are many nonequivalent calculi for a given
associative algebra.
This yields a problem of classification of first order calculi.

In our previous paper \cite{BKO} following the developments
by Pusz, Woro\-no\-wicz \cite{PW} and by Wess, Zumino \cite{WZ},
we have proposed a general algebraic formalism for a first order
calculus on an arbitrary associative algebra with a given generating
space. A basic idea was that a noncommutative differential
calculus is best
handled by means of commutation relations between generators of an
algebra and its
differentials. We assume that a bimodule of one
forms is a free right module. This allows to define partial
derivatives (vector fields).
Corresponding calculi and a differentials are said to be  {\it
coordinate
calculi} and  {\it coordinate differentials}.\\

Throughout this paper, $\F$ will denote a field.
An algebra means a linear associative unital
algebra over the field $\F$. A coordinate algebra is an algebra $R$
together with a chosen set of generators $x^1,\ld ,x^n$ (coordinates).
We assume that the coordinates are linearly independent but in
general
they do not need to be algebraically independent.
In other words, coordinate
algebra is an algebra together with its presentation.
A presentation of an algebra $R$ is an epimorphism
$\pi :\F<x^1,\ld , x^n>\rt R$, where $\F<x^1,\ld ,x^n>$ is a free
algebra generated by the variables $x^1,\ld ,x^n$.
Let $I_R=ker\pi$, then $R\cong \F<x^1,\ld , x^n>/ I_R$ where, $I_R$
is an (twosided) ideal of relations in $R$ and $\pi(x^i)$ are
generators
of $R$. By abuse of notation we continue to write $x^i$ for $\pi(x^i)$
when no confusion can arise. The free algebra is isomorphic to a
tensor
algebra $TV$ of a vector space $V={\rm lin}(x^1,\ld ,x^n)$. If $V$ is
infinite dimensional then we have infinite number of generators
$(x^1, x^2,\ld)$.

Coordinate algebra is interpreted as an algebra of (polynomial)
functions
$R=Func(X)$ on  quantum space $X$.
This resembles the situation from algebraic geometry.
Homogeneous ideals corresponds to graded algebras (projective case).
Roughly speaking coordinate algebras are quantum versions of
algebraic
varieties.

\section{Coordinate calculi and optimal algebras.}

In this section we mainly review our results from \cite{BKO}.
For the proofs we refer the reader to \cite{BKO}.

Differential or (first order differential) calculus
is a linear mapping from an algebra $R$
to a $R$- bimodule $M$ satisfying the Leibniz rule:
$$
d(uv)=d(u)v + ud(v) \eqno(1)
$$
{\bf Definition 2.1}
Let R be an algebra with a generating space $V=\*
{\rm lin}(x^1,\ld , x^n)$.
A differential $d:R\rightarrow M$ and the corresponding calculus
is said to be {\it coordinate} if the bimodule
$M$ is a free right $R$-module freely generated  by $dx^1,\ld ,
dx^n$.\\

This definition essentially depends on the generating
space $V$. The bimodule
$M$ as a right module is isomorphic to the right module $V\ot R$.
Let $v\in R$.
The left multiplication $u\mapsto vu$ is an endomorphism of the right
module $M$. Ring of all endomorphisms of any free module of rank $n$
is
isomorphic to the ring $R_{n\times n}$ of all $n$ by $n$ matrices with
entries from $R$. Therefore, we can find an algebra homomorphism
$A:R\rt R_{n\times n}$ defined by the formulae
$$
v\,dx^i=dx^k\c A(v)_k^i . \eqno(2)
$$
Clearly, it satisfies the homomorphism property
$$
A(uv)^i_k =A(u)^l_k A(v)^i_l . \eqno  (3)
$$
{\bf Example 2.2} {\it (Universal differentials)}
Let us consider the case $V$ is any supplement of the subspace
$\F~\c~1$, i.e. $R=V\oplus\F\c 1$. Of course if $R$ is not a finite
dimensional then  we have an infinite number of generators.
Nevertheless there exists a coordinate differential
with respect to this space of generators. It is exactly the universal
differential (cf. e.g. \cite{M,W}).
Let $M(R)\subset R\ot R$ denotes a kernel of the multiplication map.
Clearly, $M(R)$ is a bisubmodule. For any $u\in R$ we put
$du=1\ot u-u\ot 1$. Decomposing $u=u_0\c 1 +\sum_{i\geq 1}x^iu_i$
with $u_i\in \F$ we can find that $du=\sum_{i\geq 1}dx^iu_i$ is a
unique
decomposition. This means that $d:R\rt M(R)$ is a coordinate calculus.
The homomorphism $A$ of the
universal differential has the form
$$
A(x^i)^j_k=C^{ij}_k - \de^i_kx^j
$$
where, the scalar coefficients  $C^{ij}_k$ come from
the multiplication table for the generators
$$
x^ix^j=C^{ij}_kx^k+D^{ij}\c 1 \ .   \eqno (4)
$$
Then the homomorphism property (3) is
equivalent to  associativity constraints of the
structure constants $C^{ij}_k, D^{ij}$ vis.
$$
C^{ij}_rC^{rk}_m + D^{ij}\de^k_m =
C^{jk}_rC^{ir}_m + D^{jk}\de^i_m  .
$$\smallskip

If $d$ is a coordinate differential, then linear maps $D_k:R\rt R$
(partial
derivatives with respect to the coordinates $(x^1,\ld , x^n)$ are
uniquely
defined by the formula:
$$
d\,v=dx^k\c D_k(v) .
\eqno          (5)
$$
These maps satisfy the relations
$$
D_k(x^i)= \delta^i_k ,
\eqno     (6)
$$
and twisted derivation property
$$
D_k(uv)=D_k(u)v + A(u)_k^iD_i(v) .  \eqno    (7)
$$
By the Leibniz rule (1) we have
$$
v\,dx^i=d(vx^i) - d(v)x^i= dx^k[D_k(vx^i)-D_k(v)x^i]
$$
therefore, the partial derivatives  $D_k$ and the homomorphism $A$
are connected by the relations
$$
A(v)^i_k= D_k(vx^i)-D_k(v)x^i  \eqno (8)
$$
This shows that for a given coordinate calculus a left module
structure
on the right free module $V\ot R$ is uniquely determined by (8).
A natural question concerning an inverse problem arises here. If a
homomorphism $A$ is given, then the formula (7) allows one to
calculate
partial derivatives of a product in terms of its factors. That fact
and
formula (6) show that for a given $A$ there exists not more then one
coordinate  calculus with partial derivatives
satisfying the formulae (6) and (7). It is not clear yet whether or
not
there exists at least one differential of such a type. Thus, our first
task is to describe these homomorphisms $A$ for which there exist
coordinate differentials.\\

\noindent{\bf Theorem 2.3} \em
Let $\R=\F<x^1,\ld ,x^n>$ be a free  algebra generated by
$x^1,\ld ,x^n$ and $A^1,\ld ,A^n$ be any set of $n\times n$ matrices
over $\R$. There exists the unique coordinate differential $d_A$ of
$\R$
such that
$A(x^i)=A^i$ i.e. such that the following commutation rules are
satisfied
$$
x^idx^j=dx^k\c (A^i)^j_k .
$$\\ \em

Let now \ $R$ \ be a non-free  algebra defined by the set of
generators \ $x^1,\ld ,x^n$\ \ and an ideal $I\subset \R$ of relations
such that $R=\R / I$.\\

\noindent{\bf Definition 2.4}
An ideal $I\neq \R$ of a free algebra $\R=\F<x^1,\ld ,x^n>$ is said
to
be {\it consistent} with a homomorphism $A:\R\rt \R_{n\times n}$ if
the
the following two consistency condition are satisfied:
$$
A(I)^i_k\subseteq I \eqno (C1)
$$
$$
D_k(I)\subseteq I   \eqno (C2)
$$
where, $D_k$ are  partial derivatives for the differential $d_A$
(cf. Theorem 2.3).\\

These two conditions generalize Wess and Zumino quadratic and
linear consistency condition for quadratic algebras.\\

\noindent{\bf Theorem 2.5}\em
Let an ideal $I\neq \R$ of the free algebra $\R=\F<x^1,\ld ,x^n>$
be  consistent with a homomorphism $A:\R\rt \R_{n\times n}$. Then
the quotient algebra $R=\R / I$ has a coordinate differential with the
homomorphis $\pi (A):R\rt R_{n\times n}$, where $\pi :\R\rt R$ is
a canonical epimorphism.\em\\

\noindent{\bf Corollary 2.6}\em
Let an ideal $I\neq \R$ of the free algebra $\R=\F<x^1,\ld ,x^n>$
be generated by the set of homogeneous relations $\{ f_a\in \R, \
a\in \M\}$ of the same degree.
Assume that the homomorphism $A:\R\rt \R_{n\times n}$
acts linearly on generators i.e. $A(x^i)^j_k= \al^{ij}_{kl}x^l$.
Then the ideal $I$ is $A$--consistent or equivalently
the quotient algebra $R=\R / <f_a=0, a\in \M>$ has a coordinate
differential with the commutation rule $x^idx^j=\al^{ij}_{kl}dx^k\c
x^l$
if and only if
$$
A(f_a)^i_k=\ga^{ib}_{ka}f_b   \eqno (C1')
$$
$$
D_kf_a=0  \eqno (C2')
$$
for some scalar coefficients $\ga^{ib}_{ka}$. \\ \em

We have proved that a free algebra $\R$ admits a coordinate calculus
for
arbitrary commutation rules.
In order to define a homomorphism $A:\R\rt
\R_{n\times n}$ it is enough to set its value on generators
(cf. Theorem 2.3)
$$
(A^m)^i_k\equiv A(x^m)^i_k=\al^{mi}_k +\al^{mi}_{kj_1}x^{j_1}+
\al^{mi}_{kj_1j_2}x^{j_1}x^{j_2}+ \ld
$$
where, $\{\al^{mi}_{kj_1,\ld ,j_r}\in\F \}$ are arbitrary tensor
coefficients.
If the homomorphism $A$ preserves a degree, then it must act linearly
on
generators i.e.,
$A(x^m)^i_k=\al^{mi}_{kj}x^j$. Therefore, the homomorphism
$A$ is defined by the 2-covariant 2-contravariant tensor
$A=\al^{mi}_{kj}$. This is a homogeneous case, the most frequently
studied in the literature. In the case
$A(x^m)^i_k= \al^{mi}_k+\al^{mi}_{kj}x^j$
the partial derivatives do not increase a degree. The case
$A(x^m)^i_k = \al^{mi}_k+\de^i_kx^m$
has been considered by Dimakis and M\"uller--Hoissen \cite{DM}.

For any homomorphism $A$ there exists the largest $A$-consistent
ideal
$I(A)$ contained in the ideal $\bar{R}$ of polynomials with zero
constant
terms ($\bar{R} \oplus\F\c 1=\R$) -- the sum of all consistent ideals
of
such a type. Note that the ideal $I(A)$ need not to be the only
maximal
$A$--consistent ideal in $\R$.\\

\noindent{\bf Definition 2.7}\
The factor algebra $R(A)=\R / I(A)$ is called an {\it optimal} algebra
for a commutation rule $x^idx^j=dx^k\c A(x^i)^j_k$.\\

We shall describe the ideal $I(A)$ in the homogeneous case. \\

\noindent{\bf Theorem 2.8} \em
For any 2-covariant 2-contravariant tensor $A=\al^{ij}_{kl}$
the ideal $I(A)$ can be constructed by induction as the homogeneous
space
$I(A)=I_1(A)+I_2(A)+I_3(A)+\cdots$ in the following way:
\begin{enumerate}
   \item  $I_1(A)=0$
   \item  Assume that $I_{s-1}(A)$ has been defined and $U_s$ be a
space
   of all polynomials $m$ of degree $s$ such that $D_k(m)\in I_{s-
1}(A)$
  for all $k.\ 1\leq k\leq n$. Then $I_s(A)$ is the largest $A$-
invariant
    subspace of $U_s$.
\end{enumerate}
The ideal $I(A)$ is a maximal $A$--consistent ideal in $\R$. \em \\

We have proved that
there exists a maximal algebra which has a coordinate
differential with an arbitrary commutation rule. This is a free
algebra
$\R$. In general, there exist many algebras which have  coordinate
calculi
with a given commutation rule $A$. Each of them has the form
$\R / I$ where, $I$ is an $A$- cosistent ideal in $\R$.
The optimal algebra $R(A)$ is a minimal one in the following sense:
it does not contain $A$- consistent ideals at all. In the
homogeneous case this algebra is characterized as the unique algebra
which has no nonzero $A$-invariant subspaces with zero differentials.
We will see from the examples below that the same algebra can be
an optimal algebra for various commutation rules.


\section{Optimal versus quadratic algebras.}

In this section we shall address the existence problem for  calculi
on quadratic algebras. It turns out that this is closely related to
quadratic identities in the optimal algebra.

A coordinate algebra on $V$ is {\it quadratic} if it is
the quotient of the tensor algebra $TV$ by an ideal $<W>$ generated
by a subspace $W\subset V\ot V$.
It is customary and convenient to consider the subspace $W$ in the
form
$$
W_B={\rm lin}(x^i\ot x^j - \be^{ij}_{kl}x^k\ot x^l ) \eqno (9)
$$
with $\be^{ij}_{kl}$ being a matrix of some  endomorphism $B$ of
$V\ot V$. Elements of ${\rm End}(V\ot V)$ we shall call {\it twists}.
We will denote by $R_B$ the quadratic algebra associated with the
twist $B$
$$
R_B=TV / <W_B> \cong \F<x^1,\ld ,x^n> / <x^ix^j=\be^{ij}_{kl}x^kx^l> .
\eqno (10)
$$
This method of presentation of the algebra $R_B$ has disadvantage of
not being unique: a huge number of twists lead to the same algebra.
For example, if $\bar{W} \oplus W=V\ot V $ is the direct sum
decomposition
then the projection onto $W$ is a good candidate for $B$ (cf.
\cite{S}).
It is not our purpose to develop this point here.
An interesting case appear when the twist $B$ satisfies
the Yang--Baxter equation
$$B_{12}B_{23}B_{12}=B_{23}B_{12}B_{23} . \eqno (11)$$

In the sequel we assume
the homogeneous commutation rule $x^idx^j=\al^{ij}_{kl}dx^k\c x^l$ is
given. In this case with a homomorphis $A:\R\rt\R_{n\times n}$
one can associate a twist
determined by the same $2\times 2$ tensor
$\al^{ij}_{kl}$. By  abuse
of notation we shall denote this endomorphism by the same letter,
$A\in {\rm End}(V\ot V)$.

This yields the consistency problem between
the ideal $<W_B>$ and the commutation rule $A$.
The following proposition provides a
convenient criterion for the consistency check. \\

\noindent{\bf Proposition 3.1} {\em
Let $A,\ B$ be twists on $V$.
Then $<W_B>$ is an $A$--consistent ideal in $TV$ if and only if
the following two consistency conditions are satisfied:\\
the generalized Yang--Baxter equation
$$
A_{12}A_{23}(E-B)_{12} = (E-B)_{23}Z  \eqno (C1'')
$$
has a solution in a
$3\times 3$ tensor $Z\in {\rm End(V\ot V\ot V)}$, and
$$
(E - B)(E + A)=0  \eqno (C2'')
$$ where, $E$ is the identity twist.\\
}
{\it Proof: } The proof follows directly from the Corollary 2.6 .
\hfill $\Box$\\

Every solution of the consistency conditions (C1'') and (C2'')
leads to a coordinate calculus on the algebra $R_B$ with the
commutation
rule $A$. In particular, we have
$W_B\subseteq I_2(A)$ and $<W_B>\subseteq I(A)$ (cf. Theorem 2.8).
One can express it by
saying that the quadratic relations $x^ix^j=\be^{ij}_{kl}x^kx^l$
are fulfilled in the optimal algebra $R(A)$.
In other words,
all solutions of the consistency conditions (C1'') and (C2'') allow us
to find out all possible (quadratic) identities in the optimal algebra
$R(A)$. But in general, the quadratic algebra $R_B$ is not an optimal
algebra for the commutation rule $A$.\\

\noindent{\bf Remark 3.2}\ If the differential $d_A$ on
$R_B$ is {\it nondegenerate} in the following sense:
the kernel of $d_A$ consists only the ground field $\F$
then  $R(A)=R_B$.\\

One can also pose an inverse problem. Assume that the quadratic
algebra $R_B$ is given and we are looking for solutions $A,\ Z$ of
the
consistency
conditions. It allows us to classify all existing homogeneous
calculi on the quadratic algebra $R_B$, or equivalently to classify
all
homogeneous commutation rules for which optimal algebras satisfy the
relations $x^ix^j=\be^{ij}_{kl}x^kx^l$.\\
For example, in \cite{BKO2} we have classified all homogeneous
commutation
rules in two variables for which the optimal algebras are commutative.

Set $A, Z$ arbitrary and $B=E$. These are  trivial solutions of the
consistency conditions (C1''), (C2''). In this case $R_E\equiv TV$.
It restate in the homogenous case our previous result that
the free algebra admits a coordinate
differential with any commutation rule $A$.

Wess and Zumino have found another solutions. If twists $A, B$
satisfy the Yang--Baxter equation of the form
$$
A_{12}A_{23}B_{12} = B_{23}A_{12}A_{23} \eqno (12)
$$
then they  solve the quadratic consistency
condition (C1'').(Indeed, in this case $Z=A_{12}A_{23}$).
In particular, if $A$ satisfies the Yang--Baxter equation (11) itself
then $B=\mu A$, $\mu\in \F$ also is a  solution of (C1'').
The linear consistency condition (C2'') now becomes a minimal
polynomial
condition for $A$
$$
(E/\mu - A)(E + A) = 0 . \eqno (13)
$$
This is the Hecke equation for $A$.

More general solutions with the higher order minimal polynomial
were found by Hlavaty \cite{H}. Let $A$ be a Yang--Baxter twist
as previously. Assume $-1\in SpecA$ and $G(A)(E+A)=0$ is
some (possibly minimal) polynomial identity satisfied by $A$, where
$G(A)=\sum \mu_k A^k$ is a polynomial in the matrix $A$. Then
$B=E-G(A)$ and $Z=A_{12}A_{23}$ satisfy both consistency conditions.\\

\noindent{\bf Example 3.3}\ {\it (Manin's spaces)}\
Consider a twist $B$ in the form
$$B(x^i\ot x^j)=\be^{ji}x^j\ot x^i  \eqno (14)$$
defined by some $n$ by $n$ matrix $\be^{ij}$. In this case
$\be^{ij}_{kl}=\be^{ji}\de^i_l\de^j_k$. This choice leads to the
various Manin's spaces determined by the coefficients $\be^{ij}$
(cf. \cite{Man}).\\

\noindent{\bf Remark 3.4}\ This is an easy observation that arbitrary
three twists $A$, $B$, $C$ each of the form (14) do satisfy
the Yang--Baxter equation $A_{12}C_{23}B_{12}=B_{23}C_{12}A_{23}$.\\

Let us consider the commutation rule $A$ of the form (14), i.e.:
$x^jdx^i=dx^i\c\al^{ij}x^j$, with some matrix $\al^{ij}$.
Look for
quadratic relations given by the formula (14)
in the optimal algebra.
Due to the Remark 3.4 the quadratic consistency condition (C1'')
is automatically satisfied in the form (12).
Substituting to the linear
condition (C2'') one gets
$$
(\al^{ji}-\be^{ji})\de^j_m\de^i_p + (1-
\be^{ji}\al^{ij})\de^i_m\de^j_p=0
$$
Consider two cases. For $i=j$ we have the condition
$$(\al^{ii}+1)(1-\be^{ii})=0 .$$
It means that $\be^{ii}$ has to be $1$ whenever $\al^{ii}\neq -1$.\\
For $i\neq j$ we obtain two conditions
$$\al^{ij}=\be^{ij} \ \ \mbox{and} \ \ \ \al^{ij}\al^{ji}=1 .$$
Thus we get the following result.
If the matrix $\al^{ij}$ does satisfy the property $\al^{ij}\al^{ji}
=1$ for $i\neq j$ then $\be^{ij}=\al^{ij}$ for $i\neq j$ and
$\be^{ii}=1$
solves both consistency conditions.
In this case  the relations $\al^{ij}x^ix^j=x^jx^i$ for $i\neq j$
are satisfied in the optimal algebra $R(A)$. It means that the Manin
space
$$
R_B=\F<x^1,\ld ,x^n>/ <\al^{ij}x^ix^j=x^jx^i, \ i<j > \eqno (15)
$$
have a coordinate calculus with our commutation rule $A$.
For this algebra to be optimal by the Remark 3.2
it is enough to check that this differential is nondegenerate.
It was done in \cite{BKO2}. Summarizing, we have\\

\noindent{\bf Example 3.5}\ {\it (Calculi on Manin's spaces)}\
Consider the diagonal commutation rule $x^jdx^i=\al^{ij}dx^i\c x^j$
with $\al^{ij}\al^{ji}=1$ , for $i\neq j$. There are two cases.
\begin{itemize}
\item[(A)]
If
none of the coefficients $q^{ii}$ is a root of a polynomial of the
type
$\la^{[m]}\doteq \la^{m-1}+\la^{m-2}+\cdots +1$
then the optimal algebra $R(A)$ for this commutation rule
has the form (15).
\item[(B)]
If $(\al^{ii})^{[m_i]}=0,\ 1\leq i\leq s$ with the minimal $m_i$ then
\\
$R(A)=R_B/ <(x^i)^{m_i}=0,\ 1\leq i\leq s>=$. \\
$$\F<x^1,\ld ,x^n>/<\al^{ij}x^ix^j=x^jx^i,\ i<j,\ (x^i)^{m_i}=0,\
1\leq i\leq s > \eqno (16) $$
\end{itemize}
It should be noted that in the case (B) the optimal algebra is not
in general a quadratic algebra. The paragrassmann case $m_i=m>2$
has been studied in \cite{FIK}.\\

Observe that in the case (A), the diagonal elements $\al^{ii}$ have
no influence on the form (15) of the algebra $R(A)$. But different
$\al^{ii}$ lead to different commutation rule and different
differentials.
This is the easiest way to see that the same algebra can be an optimal
algebra for different calculi.

It is rather expected result that
the polynomial algebra (free commutative algebra)
is an optimal algebra for the Newton--Leibniz calculus $(\al^{ij}=1)$,
while the Grassmann algebra is optimal for the supersymmetric
calculus ($\al^{ij}=-1$).\\

\section{Higher order calculi.}

\noindent{\bf Definition 4.1.}\ A {graded differential algebra} (GDA
for
short)
is a $\N$--graded algebra $$\Om=\oplus_{m\in \N}\Om^m$$
and a homogeneous $\F$--linear mapping $d$ of degree $1$
$$d:\Om\rt \Om, \ \ d=\oplus_{m\in\N}d_m, \ \ d_m:\Om^m\rt \Om^{m+1}$$
such that:
$$d^2=0$$
$$d(\om_1\c \om_2)=d\om_1\c \om_2+(-1)^{deg\om_1}\om_1\c \om_2$$
for each $\om_1, \ \om_2 \in \Om$, $\om_1$ homogeneous.\\\smallskip

In particular, it follows from the above definition that $\Om^0$ is an
algebra, $\Om$ is an $\Om^0$-- algebra, each $\Om^n$ is an
$\Om^0$-- bimodule and $d_0:\Om^0\rt \Om^1$ is a (first order)
differential. This fact can be expressed by saying that the GDA $\Om$
is an extension of the first order calculus $d_0$ to higher order.
Since GDA is also a complex ($d^2=0$) one can think of it as a
noncommutative or quantum generalization of a de Rham complex.
Conversely, if the first order differential $d:R\rt M$ is given,
there exist in general many extension of it to a GDA $\Om$ in such
a way that $\Om^0=R$, $\Om^1=M$ and $d_0=d$. Several such
constructions
have been proposed in the literature.\\

\noindent{\bf Example 4.2.}\ Wornowicz's external algebra formalism
for quantum groups \cite{W}.\\

\noindent{\bf Example 4.3.}\ {\it (Universal GDAs)}\
Let $d:R\rt M$ be a first order differential
such that the bimodule $M$ is generated by differentials i.e.
$M=dR\c R$. Maltsiniotis \cite{M} has shown that there exists a GDA
$\Om(R, M)$ which satisfies certain universal property:
each GDA extension of $d:R\rt M=dR\c R$ can be obtained from $\Om(R,
M)$
via a standard quotient construction. In particular, if $d:R\rt M(R)$
is
a universal first order differential (cf. Example 2.2) then the
universal GDA $\Om(R)\equiv \Om(R, M(R))$ is known as a differential
envelope of the algebra $R$.\\

Here, we propose a general construction
of GDA $\Om_A(R, M)$
for a coordinate differential $d:R\rt M$ with a commutation rule
$x^idx^j=dx^k\c A(x^i)^j_k$.
This construction is related to an external algebra formalism.
The detailed discussion and proof will appear
in a forthcoming publication \cite{BK}.\\

Let $T_RM=R\oplus M\oplus M\ot_RM\oplus\ld$ be a tensor algebra of
the
bimodule $M$. $T_RM$ is a free right module freely generated by
elements
of the form $dx^{i_1}\ot\ld\ot dx^{i_m}$, \ $m\in\N$, i.e.
$T_RM=T(dV)\ot_\F R$
where, $T(dV)$ is a tensor algebra of the vector space
$dV={\rm lin}(dx^1,\ld ,dx^n)$. Define a linear mapping $d:T_RM\rt
T_RM$ by
$$
d(dx^{i_1}\ot\ld\ot dx^{i_m}\c v_{i_1\ld i_m})=
(-1)^mdx^{i_1}\ot\ld\ot dx^{i_m}\ot dv_{i_1\ld i_m}       \eqno (17)
$$
Note that $(T_RM, d)$ is not a GDA.\\

\noindent{\bf Theorem 4.4.}\ {\em
Let $J_A$ denotes an $R$- homogeneous ideal in the $R$- algebra $T_RM$
generated by the $n^2$ elements
$$dx^i\ot dx^j+dx^k\ot dx^l\c D_l(A(x^i)^j_k)$$
Let $\Om_A(R, M)=T_RM / J_A$. Then
\begin{description}
\item[1.] (GDA)\\
$$\Om_R(R, M)=R\oplus M\oplus\Om^2_A(R, M)\oplus\ld$$
is a GDA where, now $d$ is a mapping induced from (17).
\item[2.] (existence of wedge product)\\
In the homogeneous case $A(x^i)^j_k=\al^{ij}_{kl}x^l$ one has
$$\Om^m_A(R, M)=\La^m_A(dV)\c R$$
i.e. $\Om^m_A(R, M)$ is generated by elements from $\La^m_A(dV)$
where,
$$\La_A(dV)=\oplus_{m\in\N}\La^m_A(dV)=T(dV) /
<dx^i\ot dx^j+dx^k\ot dx^l\al^{ji}_{kl}>$$
\item[3.] (external algebra structure)\\
Moreover, the equality
$$\Om^m_A(R, M)=\La^m_A(dV)\ot_\F R$$
holds, i.e. the bimodule $\Om^m_A(R, M)$ is a free right module
if and only if the following generalized Yang--Baxter equation has a
solution in a $3\times 3$ tensor $Y$:
$$A_{12}A_{23}A_{12} - A_{23}A_{12}A_{23}=(E+A)_{12})Y$$
here, $A=\al^{ij}_{kl}$ .

\end{description}
}\medskip

\noindent{\bf Acknowledgments. }\
A part of this work is based on joint papers and a fruitful
collaboration
with Z. Oziewicz. It is a great pleasure to thank him.

\end{document}